\documentclass[11pt,twoside]{article}

\usepackage{amsmath}
\usepackage[all]{xy}

\newtheorem{thm}{Theorem}[section]

  \newtheorem{defn}[thm]{Definition}

\begin{document}


 \centerline{\large{\textbf{ The fiber homotopically equivalent      }}}
  \centerline{\large{\textbf{relation for fiber bundles  }}}

  \centerline{{\emph{\textbf{Suliman.D and Adem K{\i}l{\i}\c cman}}}}



\begin{abstract}
\footnotesize{   In this paper, we will give the polyhedron
property  role to satisfy fiber homotopically equivalent relation
in fiber bundle theory over suspensions
of polyhedron space. \\
\textbf{\emph{Keywords}}:  fibration; fiber bundle; polyhedron;   homotopy.\\
 ------------------------------------\\
 \emph{AMS classification}: Primary 53D28, 55R10 , 52B70, 14F35}

\end{abstract}

\section{Introduction}

 Throughout in this paper the word "space" means Hausdorff
space, the word "H-fibration" means an onto  regular Hurewicz
fibration, and the space of parameterized paths in any space with
the compact-open topology. Also we mean by symbol $\simeq$
"homotopy relation". For a path $\alpha$ in any space $X$, we
denote the inverse path of $\alpha$ by $\overline{\alpha}$.  In
  H-fibration   $P:E \longrightarrow B$, we denote by $P|A$ the
H-fibration $P|_{P^{-1}(A)}: P^{-1}(A)\longrightarrow A$ defined
by $ P|_{P^{-1}(A)}(e)=P(e)$ for $e\in P^{-1}(A)$, where $A$ is a
subspace of $B$.

\begin{thm}\label{rsul}
\emph{ \cite{Spa}   Let $X$, $Y$ and $Z$   be topological spaces.
If $X$ is locally compact and regular space then the map
$H:Z\longrightarrow Y^X$
 always gives rise to the map $F:Z\times X\longrightarrow Y$ by defining
 $F(z,x)=H(z)(x)$ for all $x\in X, z\in Z$.  }
    \end{thm}

Amin and Alinor \cite{AmA2}  introduced the notion of the
$Sf-$function in the theory for fibrations as follows:

\begin{defn}\label{d:521}
\emph{Let $P:E\longrightarrow B$ be a H-fibration with the a
lifting function $\lambda$ and fiber spaces $F_{b_o}=P^{-1}(b_o)$,
where $b_o\in B$.  The \emph{$Sf-$function for   $P$ induced by
$\lambda$} is a  map $\Theta_{\lambda}:L(B,b_o)\times
F_{b_o}\longrightarrow F_{b_o}$ defined by
$\Theta_{\lambda}(\alpha,e)=\lambda(e,\alpha)(1)$ for all $e\in
F_{b_o}, \alpha\in L(B,b_o)$, where $L(B,b_o)$ is the set of all
loops in $B$ based $b_o$.}
\end{defn}

\begin{defn}\label{d:553}
\emph{Let $P_{1}:E_{1}\longrightarrow B$ and
$P_{2}:E_{2}\longrightarrow B$ be two H-fibrations with fiber
spaces $F^{1}_{b_o}=P^{-1}_{1}(b_o)$ and
$F^{2}_{b_o}=P^{-1}_{2}(b_o)$, where $b_o\in B$. The
$Sf-$functions $\Theta_{\lambda_1}:L(B,b_o)\times
F^{1}_{b_o}\longrightarrow F^{1}_{b_o}$ and
$\Theta_{\lambda_2}:L(B,b_o)\times F^{2}_{b_o}\longrightarrow
F^{2}_{b_o}$ are said to be \emph{conjugate }if there is $g\in
H(F^{1}_{b_o},F^{2}_{b_o})$ such that $\Theta_{\lambda_1}\simeq
\overline{g}\circ \Theta_{\lambda_2}\circ(id_{L(B,b_o)}\times g)$,
where
 $H(F^{1}_{b_o},F^{2}_{b_o})$ is the set of all homotopy equivalences from $ F^{1}_{b_o}$ into $F^{2}_{b_o}$
 and $\overline{g}$ denotes  the homotopy inverse of $g$. }
\end{defn}

\begin{thm}\label{t:559}
\emph{\cite{suliman} Let $ B $ be a polyhedron and be the union of
two subpolyhedrons  $B_1$ and $B_2$ such that $B_1$ is a
contractible in $B$   to a point $b_o\in B_{1}\cap B_{2}$ leaves
$b_o$ fixed and $B_2$ is also a contractible to $b_o$ and
$B_{1}\cap B_{2}$ be subpolyhedron of $B$. If  $P_{1}: E_{1}
\longrightarrow
 B$  and
$P_{2}:E_{2} \longrightarrow B $ are two H-fibrations with
conjugate $Sf-$functions $\Theta_{\lambda_1}$ and
$\Theta_{\lambda_2}$by $g\in H(F^{1}_{b_o},F^{2}_{b_o})$, then
$P_{1}$ and $P_{2}$ are fiber homotopy equivalent.}
\end{thm}


\begin{defn}\label{fiber1}
\emph{Let $E$, $B$, and $F$ be spaces. Let $P:E\longrightarrow B$
be a map of $E$ onto $B$ and $G$ be group of all homeomorphisms of
$F$ onto $F$ with as a binary usual composition operation $\circ$.
Then $\gamma=(E,P,B,F,G)$ is said to be \emph{ a fiber bundle over
a base $B$} if there is an open covering  $\{V_j:j\in \wedge\}$ of
$B$, where $\wedge$ is index set and for each $j\in \wedge$, there
is a homeomorphism $ \theta_j:V_j \times F\longrightarrow
P^{-1}(V_j)$ ( called coordinate function ) such that:\newline
    1- $P[\theta_j(b,y)]=b$ for all $b\in V_j, y\in
    F$.\newline
    2- For each pair $i,j\in \wedge$ and $b\in V_i\cap V_j$, the homeomorphism
    $\theta^{-1}_{jb}\circ \theta_{ib}:F\longrightarrow F$ corresponds to an element of $G$, where
     $\theta_{kb}: F\longrightarrow P^{-1}(b)$ defined by
     $\theta_{kb}(y)=\theta_k(b,y)$  for all $b\in V_k, y\in
    F$ and $k=i,j$.\newline
    3-  For each pair $i,j\in \wedge$, the function
    $g_{ij}:V_i\cap V_j\longrightarrow G$ given by $g_{ij}(b)=\theta^{-1}_{jb}\circ
    \theta_{ib}$ is continuous ( called coordinate transformation).}
\end{defn}

 In fiber bundle $\gamma=(E,P,B,F,G)$, we shall denote the identity element of a group $G$ by $e$,
        the inverse element $g\in G$ by  $g^{-1}$,  we mean by a map $k:(X,x_o)\longrightarrow (Y,y_o)$ a map $k$ of $X$ into $Y$ and
        $k(x_o)=y_o$.

It's clear that every fiber bundle $\gamma=(E,P,B,F,G)$ is a
regular local Hurewicz fibration and if
      $B$ is paracompact then $\gamma$ is  Hurewicz fibration.
\begin{thm}\label{fiber4}
 \emph{\cite{James} Let $S^n$ be the $n-$sphere in $R^{n+1}$, where $n>0$ is a positive integer. Then
 for a fiber bundle $\gamma=(E,P,S^n,F,G)$, there is a characteristic map $\mu:(S^{n-1}, x_o)\longrightarrow (G,e)$.}
 \end{thm}

Now we recall the Dold's theorem in fiber bundles over sphere
$S^n$ as follows:
\begin{thm}\label{fiber5}
\emph{[The Dold's theorem] \newline Let $\gamma =(E ,P ,S^n,F ,G )
$ and $\gamma'=(E',P',S^n,F',G')$  be two fiber bundles over
sphere $S^n$ with locally compact fibers $F $and $F'$. Let $\mu
:(S^{n-1}, x_o)\longrightarrow (G ,e )$  and $\mu':(S^{n-1},
x_o)\longrightarrow (G',e')$ be characteristic maps of $\gamma $
and $\gamma'$, respectively and let $i :G \longrightarrow H(F , F
)$ and $i':G'\longrightarrow H(F', F')$ be the inclusion maps.
Then $\gamma $ and $\gamma'$ are fiber homotopy equivalent if and
only if there is homotopy equivalence $g: F \longrightarrow F'$
such that the maps \[ q (x)=  g\circ (i \circ \mu )(x)\circ
\overline{g}\quad \mbox{and}\quad q'(x)=(i'\circ \mu')(x)\] of
$S^{n-1}$ into $H(F',F')$ are homotopic.}
\end{thm}

\begin{defn}\label{fiber6}
\emph{Let $X$ be any space and $x_o$ be a base point in $X$. The
\emph{suspension } $S(X)$  of $X$ is defined to be the quotient
space of $X\times I$ in which for all $x\in X$, $(x,0)$ identified
to $(x_o,0)$,  $(x,1)$ identified to $(x_o,1)$, and  $X\times
\{1/2\}$ identified to $X$. }
\end{defn}
 The Dold's theorem remains valid if we use suspensions of
polyhedron spaces  instead of $n-$spheres $S^n$ for the
      base of bundles,  \cite{Husem}.

\begin{thm}\label{fiber8}
 \emph{ \cite{James} The suspension $S(S^n)$ of $n-$sphere $S^n$ is homeomorphic to sphere $S^{n+1}$, where $n>0$ is a positive integer.}
      \end{thm}
\begin{thm}\label{fiber9}
 \emph{  \cite{Husem}  If $X$ is a polyhedron space, then $S(X)$ is also  polyhedron space.}
      \end{thm}

\section{\emph{Sf}-function  and fiber bundles}
In this paper, we will show the \emph{Sf}-function role to satisfy
fiber homotopically equivalent relation in fiber bundle theory, in
particular,  we will the equivalently between our
theorem(\ref{t:559})  and Dold's theorem over suspensions of
polyhedron spaces.
\newline

As mentioned perviously that   Dold's theorem remains valid if we
use suspensions of polyhedron space  instead of $n-$spheres $S^n$
for the
      base of bundles. Hence theorem  (\ref{fiber4}) also  remains valid  with suspensions of polyhedron spaces.
      That is,
 for a fiber bundle $\gamma=(E,f,S(X),F,G)$ over suspension $S(X)$ of polyhedron space $X$,
 there is a characteristic map $\mu:(X, x_o)\longrightarrow
 (G,e)$.  In the following theorem, we will  prove that the converse of this
 theorem is also true for suspension $S(X)$ of any space$X$.
\begin{thm}\label{t642}
\emph{Let $G$ be a group of all homeomorphisms of space $F$ with
as binary usual composition operation $\circ$ and $X$ be
 any space. If there is a map $\mu:(X, x_o)\longrightarrow
(G,e)$, then there is bundle $E$ over  $S(X)$ and a map
$P:E\longrightarrow S(X)$ such that $\gamma=(E,P,S(X),F,G)$ is
fiber bundle.}
\end{thm}
\textbf{Proof.} By the definition of $S(X)$,  we can consider
$S(X)$ as the union of two cones, one of them is $C_0(X)=X\times
[0,1/2]$ with $(x,0)$ identified to $(x_o,0)$ and
       the other  $C_1(X)=X\times [ 1/2,1]$ with $(x,1)$ identified to
       $(x_o,1)$. Let $V_1=C_0(X)$ and $V_2=C_1(X)$, then $S(X)= V_1\cup V_2$ and there is
        a retraction $r:V_1\cap V_2\longrightarrow X$.
 Now define maps
 \[g_{ii}:V_i\longrightarrow G, \quad g_{ii}(x)=e\quad \forall
 \mbox{ }x\in V_i, (i=1,2),\]
 \[g_{12}:V_1\cap V_2\longrightarrow G, \quad g_{12}(x)=(\mu \circ r)(x)\quad \forall
 \mbox{ }x\in V_1\cap V_2,\]and
 \[g_{21}:V_1\cap V_2\longrightarrow G, \quad g_{21}(x)=[g_{12}(x)]^{-1}\quad \forall
 \mbox{ }x\in V_1\cap V_2.\]
 Let $J=\{1,2\}$ be a space with the discrete topology and let
 $T\subset S(X) \times F\times J$ be the set defined by
 \[T=\{(x,y,j):x\in V_j, y\in F, j\in J\}.\]
 Define an equivalent relation $\equiv$ on $T$ by
 \[(x_1,y_1,j)\equiv (x_2,y_2,k) \Longleftrightarrow x_1=x_2 \quad
 \mbox{and} \quad g_{kj}(x_1)(y_1)=y_2,\] where $(x_1,y_1,j),(x_2,y_2,k)
 \in T$. Then put $E$ be the set of equivalence classes so
 obtained with the quotient topology.
 Hence define a map $P:E\longrightarrow S(X)$ by
 \[ P([(x,y,j)])=x \quad \forall
 \mbox{ }[(x,y,j)]\in E,\] and the maps $
\theta_j:V_j \times F\longrightarrow P^{-1}(V_j)$ defined by
\[\theta_j(x,y)=[(x,y,j)]\quad \forall
 \mbox{ } (x,y)\in V_j \times F.\] Hence it's clear that $\gamma=(E,P,S(X),F,G)$ is fiber
 bundle.\quad $ \diamond$
\begin{thm}\label{t643}
\emph{Let $P:E\longrightarrow B$ be H-fibration with locally
compact fiber $F=P^{-1}(b_o)$, where $b_o\in B$. Then the function
$\phi:L(B,b_o)\longrightarrow F^F$ given by \[
\phi(w)(x)=\Theta_{\lambda}(w,x)\quad \forall \mbox{ } w\in
L(B,b_o), x\in F,\] is a continuous function of  $L(B,b_o)$ into
$H(F,F)$.}
\end{thm}
\textbf{Proof.} Since $F$ is a locally compact and Hausdorff, then
$F$ is regular. Hence  by theorem (\ref{rsul}),  the function
$\phi$ is continuous. Now we will prove that for $w\in L(B,b_o)$,
$\phi(w)$ is homotopy equivalence from $F$ into $F$. For $w\in
L(B,b_o)$, we can define a map $ \overline{\phi(w})
:F\longrightarrow F$ by
\[
 \overline{\phi(w)} (x)=\Theta_{\lambda}(\overline{w},x)\quad \forall \mbox{
}   x\in F.\] Then we get that \[[\phi(w)\circ
 \overline{\phi(w)} ](x)=\lambda[\lambda(x,\overline{w})(1),w](1)\quad
\forall \mbox{ } x\in F,\] and
\[ \overline{\phi(w)} \circ \phi(w)](x)=\lambda[\lambda(x,w)(1),\overline{w}](1)\quad
\forall \mbox{ } x\in F.\] Then by lemma (2.2) in \cite{suliman},
\[\phi(w)\circ  \overline{\phi(w)} \simeq i_F \quad \mbox{and}\quad
 \overline{\phi(w)} \circ \phi(w)\simeq i_F.\] Hence
$\phi(w):F\longrightarrow F$ is homotopy equivalence, that is,
$\phi(w)\in H(F,F)$. Therefore $\phi$ is a map from $L(B,b_o)$
into $H(F,F)$.\quad $\diamond$
\begin{thm}\label{t644}
\emph{Let $\gamma=(E,P,B,F,G)$ be a fiber bundle admits a lifting
function $ \lambda$ with locally compact fiber $F$. Then the
function $\phi:L(B,b_o)\longrightarrow F^F$ given by \[
\phi(w)(x)=\Theta_{\lambda}(w,x)\quad \forall \mbox{ } w\in
L(B,b_o), x\in F,\] is a map from $L(B,b_o)$ into $G$.}
\end{thm}
\textbf{Proof.} Since $F$ is a locally compact, then by theorem
(\ref{rsul}),  the function $\phi$ is continuous.  For $w\in B^I$,
let $F_{w(0)}=P^{-1}(w(0))$ and let $F_{w(1)}=P^{-1}(w(1))$. Then
the map $A: F_{w(0)}\longrightarrow F_{w(1)}$ given by
\[A(x)=\lambda(x,w)(1)\quad \forall \mbox{ }x\in  F_{w(0)},\] is homeomorphism since it is obtained
from the compositions of coordinate functions which are
homeomorphisms. Hence  $\phi$ is a map from $L(B,b_o)$ into $G$.
\quad $\diamond$
\newline

To prove the equivalently between   theorem  (\ref{t:559}) and
  Dold's theorem, we first have to  rephrase  theorem
(\ref{t:559}) for two H-fibrations over a common suspension base
as follow: \newline

For any space $E$ with fixed point $x_o\in E$, there is a
\emph{conical map} $\psi:E\longrightarrow L(S(E),x_o)$ gives
       as follows:\newline
       consider $S(E)$ as the union of two cones one of them $C_0(E)=E\times [0,1/2]$ with $(x,0)$ identified to $(x_o,0)$ and
       the other  $C_1(E)=E\times [ 1/2,1]$ with $(x,1)$ identified to $(x_o,1)$. $C_0(E)$ can be contracted on itself to $x_o$ leaving
       $x_o$ fixed and similarly for $C_1(E)$. And for $x\in E$, let $w_0(x)$ be path between $x$ and $x_o$ in $C_0(E)$
       and let $w_1(x)$ be path between $x$ and $x_o$ in $C_1(E)$. Define the \emph{conical map }$\psi:E\longrightarrow L(S(E),x_o)$ by
       \[\psi(x)=\overline{w_1(x)}  \star w_0(x)\quad \forall \mbox{ }x\in  E.\]
\newline

Let $\gamma =(E ,P ,S(X),F ,G ) $  and $
\gamma'=(E',P',S(X),F',G')$  be two H-fibrations over a common
suspension base $S(X)$  of a polyhedron space $X$ with  with
locally compact fibers $F $and $F'$. In Figure 6.1, let \[\mu :(X,
x_o)\longrightarrow (G ,e )\quad \mbox{and}\quad\mu':(X,
x_o)\longrightarrow (G',e')\] be characteristic maps of $\gamma $
and $\gamma'$, respectively. And
\[i :G \longrightarrow H(F , F )\quad \mbox{and}\quad
i':G'\longrightarrow H(F', F')\] be the inclusion maps. From
theorems (\ref{t643}) and (\ref{t644}), then
  theorem  (\ref{t:559}) and   Dold's theorem can now be
compared.

\begin{displaymath} \xymatrix{
 &&&& L(S(X),x_o) \ar[dllll]_{\phi'} \ar[drrrr]^{\phi}  &\\
 G'  \ar[d]_{i'}  &&&& X\ar[u]^{\psi} \ar[llll]^{\mu'} \ar[rrrr]_{\mu}   &&&& G \ar[d]_{i}& \\
H(F',F')                              &&&&&&&& H(F,F)
\ar[llllllll]^{T(f)=g\circ f\circ [\overline{g}] \quad \forall
\mbox{ } f\in H(F,F)} }
\end{displaymath}
 Where $g\in H(F,F')$ and $\psi$ is the conical map.\newline

 Hence  theorem  (\ref{t:559})  can then be restated in terms of $\phi$,
 $\phi'$, and $\psi$ as follows:\newline
 Two H-fibrations  $\gamma =(E ,P ,S(X),F ,G ) $ and $
\gamma'=(E',P',S(X),F',G')$ are  fiber homotopy equivalent if and
only if there is $g\in H(F,F')$ such that two maps
\[m(x)=g\circ i\circ \phi[\psi(x)]\circ \overline{g}\quad \forall \mbox{
} x\in X,\] and
\[m'(x)= i'\circ \phi'[\psi(x)]\quad \forall \mbox{
} x\in X,\] from $X$ into $H(F',F')$ are homotopic.\newline
\newline

  Now if  $\phi \circ \psi\simeq \mu$ and  $\phi'
\circ \psi\simeq \mu'$, then theorem  (\ref{t:559}) and  Dold's
theorem are equivalent. In the following  theorem, we will prove
the desired.

\begin{thm}\label{t647}
\emph{Let $\gamma=(E,P,S(X),F,G)$ be a fiber bundle over
suspension $S(X)$ of a polyhedron space $X$   with locally compact
fiber $F$ and admits a lifting function $ \lambda$. And let
$\phi:L(S(X).x_o)\longrightarrow G$ be a map given by
\[ \phi(\beta)(x)=\Theta_{\lambda}(\beta,x) \quad \forall \mbox{ } \beta\in
L(S(X),x_o), x\in F.\] Then $\phi \circ \psi\simeq \mu$, where
$\mu:(X,x_o)\longrightarrow (G,e)$ is the characteristic map of
$\gamma$ and $\psi$ is the conical map.}
\end{thm}
\textbf{Proof.} Let $B=S(X)=C_0(X)\cup C_1(X)$, $B_1=C_0(X)$, and
$B_2=C_1(X)$. It's clear that $x_o\in X=B_1\cap B_2$. Now define
maps
 \[g_{ii}:B_i\longrightarrow G, \quad g_{ii}(x)=e\quad \forall
 \mbox{ }x\in B_i, (i=1,2),\]
 \[g_{12}:X\longrightarrow G, \quad g_{12}(x)= \mu(x)\quad \forall
 \mbox{ }x\in X,\]and
 \[g_{21}:X\longrightarrow G, \quad g_{21}(x)=[\mu(x)]^{-1}\quad \forall
 \mbox{ }x\in X.\]
 Let $J=\{1,2\}$ be a space with the discrete topology and let
 $T\subset S(X) \times F\times J$ be the set defined by
 \[T=\{(x,y,j):x\in B_j, y\in F, j\in J\}.\]
 Define an equivalent relation $\equiv$ on $T$ by
 \[(x_1,y_1,j)\equiv (x_2,y_2,k) \Longleftrightarrow x_1=x_2 \quad
 \mbox{and} \quad g_{kj}(x_1)(y_1)=y_2,\] where $(x_1,y_1,j),(x_2,y_2,k)
 \in T$.\newline

Recall theorem (\ref{t642}) that points of $E$ are identified to
the equivalent classes of all trips $(x,y,j)\in T$. Hence for
$j=1,2$, the maps $\epsilon_j:B_j\times F\longrightarrow
P^{-1}(B_j)$ given by
\begin{eqnarray}\label{suliman1}
  \epsilon_j(x,y)&=&[(x,y,j)]\quad \forall
 \mbox{ } (x,y)\in B_j \times F,
\end{eqnarray}
are fiber homeomorphisms which
  denote  the equivalence class
 of the triple $(x,y,j)$. Hence put $y=[(x_o,y,j)]$, where
 $j=1,2$.\newline

Hence we can define lifting functions $\lambda_1$ and $\lambda_2$
for fibrations $P|B_1$ and $P|B_2$, respectively,  as
follow:\newline
\[\bigtriangleup^1 P=\{(e,w)\in E\times B^I_1 :
P(e)=w(0)\},\]
 \[\bigtriangleup^2 P=\{(e,w)\in E\times B^I_2 :
P(e)=w(0)\},\]
\[\overline{\bigtriangleup}  P=\{(\beta,x)\in L(B,b_o)\times F :
\beta=w_2\star w_1,\]
\[ \quad\mbox{where}\mbox{ }w_i\in B^I_i
(i=1,2)\quad \mbox{and}\quad \beta(1/2)=w_2(1)=w_1(0)\in B_1\cap
B_2\},\]
\begin{eqnarray}\label{suliman2}
 \lambda_1(e,w)&=&\epsilon_1[w(t),(\pi_2\circ \epsilon^{-1}_1)(e)]\quad \forall  (e,w)\in \bigtriangleup^1 P,
 \end{eqnarray}
and
\begin{eqnarray}\label{suliman3}
\lambda_2(e,w)&=&\epsilon_2[w(t),(\pi_2\circ
\epsilon^{-1}_2)(e)]\quad \forall  (e,w)\in \bigtriangleup^2 P,
\end{eqnarray}
where $\pi_2$ is the second projection.\newline

Since $\lambda$ is  lifting function for $\gamma$, then it is also
 lifting function for $P|B_1$ and
$P|B_2$. Hence $\lambda \simeq \lambda_1$ on  $\bigtriangleup^1 P$
and $\lambda \simeq \lambda_2$ on  $\bigtriangleup^2 P$.  Hence
the map $\overline{\phi}:\overline{\bigtriangleup}
P\longrightarrow F$ given by
\begin{eqnarray}\label{suliman4}
  \overline{\phi}(\beta,x)&=&\lambda_1[\lambda_2(x,w_2)(1),w_1](1)
\quad\forall\mbox{ } (\beta,x)\in \overline{\bigtriangleup} P,
\end{eqnarray}
is homotopic to the map
$\widetilde{\phi}:\overline{\bigtriangleup} P\longrightarrow F$
defined
 by
 \begin{eqnarray}\label{suliman5}
  \widetilde{\phi}(\beta,x)&=&\lambda [\lambda (x,w_2)(1),w_1](1)
\quad\forall\mbox{ } (\beta,x)\in \overline{\bigtriangleup} P.
\end{eqnarray}
That is,  $ \widetilde{\phi}\simeq \overline{\phi}$ and by lemma
(2.2) in \cite{suliman},  $ \widetilde{\phi}\simeq
\Theta_{\lambda}$. Hence
\begin{eqnarray}\label{suliman6}
 \Theta_{\lambda}&\simeq&
\overline{\phi}.
\end{eqnarray}

 Now from the equations 1,2,3, we have that for $\beta\in L(S(X),x_o)$,
 \[\beta=w_2 \star w_1, \quad w_1\in B^I_1, w_2\in B^I_2,\]
 \[\lambda_1(y,w_1)(1)=[(w_1(1),y,1)],\] and
 \begin{eqnarray*}
   \lambda_2(y,w_2)(1)&=&[(w_2(1),y,2)]\\
     &=& [(w_2(1),\mu(w_2(1))(y),1)].
 \end{eqnarray*}
Hence
\begin{eqnarray*}
  \overline{\phi}(\beta,x)&=&\lambda_1[\lambda_2(x,w_2)(1),w_1](1) \\
   &=& [(x_o,\mu(w_2(1))(y),1)] \\
    &=& \mu(w_2(1))(y) \\
    &=& \mu(\beta(1/2))(y).
  \end{eqnarray*}
  Let $\overline{L}(S(X),x_o)$ be the projection of
  $\overline{\bigtriangleup}P$ on $L(S(X),x_o)$ and
  $\overline{\overline{\phi}}:\overline{L}(S(X),x_o)\longrightarrow G$ be a map given by
\[  \overline{\overline{\phi}}(\beta)(x)=\overline{\phi}(\beta,x)\quad\forall\mbox{ }  \beta \in \overline{L}(S(X),x_o), x\in F.\]
Then by the equation 6, $\phi \simeq \overline{\overline{\phi}}$
and
\[(\overline{\overline{\phi}}\circ \psi)(x)=\overline{\overline{\phi}}[\psi(x)]=\mu[\psi(x)(1/2)]=\mu(x),\]
that is, $\overline{\overline{\phi}}\circ \psi=\mu$. Hence
$\phi\circ \psi\simeq\mu$. Therefore theorem  (\ref{t:559}) and
the Dold's theorem are equivalent.
 \quad $\diamond$

\bibliographystyle{plain}


\end{document}